\documentclass[12pt]{amsart}

\usepackage{amscd,amsmath,amsthm,amssymb}
\usepackage{color}
\usepackage{tikz}
\usepackage{url}
\usepackage{circuitikz}
\usepackage{float}
\usepackage{lipsum}

\newtheorem{Theorem}{Theorem}[section]
\newtheorem{Lemma}[Theorem]{Lemma}
\newtheorem{Corollary}[Theorem]{Corollary}
\newtheorem{Proposition}[Theorem]{Proposition}
\newtheorem{Remark}[Theorem]{Remark}

\newtheorem{Example}[Theorem]{Example}

\newtheorem{Question}[Theorem]{Question}

\def\qed{\ifhmode\textqed\fi
	\ifmmode\ifinner\quad\qedsymbol\else\dispqed\fi\fi}
\def\textqed{\unskip\nobreak\penalty50
	\hskip2em\hbox{}\nobreak\hfill\qedsymbol
	\parfillskip=0pt \finalhyphendemerits=0}
\def\dispqed{\rlap{\qquad\qedsymbol}}

\def\m{\mathfrak{m}}

\def\pd{\textup{proj.dim}}

\begin{document}
	
	\title{Minimal cellular resolutions of complementary edge ideals}
	\author{Somayeh Moradi}

	\address{Somayeh Moradi, Department of Mathematics, Faculty of Science, Ilam University, P.O.Box 69315-516, Ilam, Iran}
	\email{so.moradi@ilam.ac.ir}
	
	
\begin{abstract}
We construct an explicit minimal cellular  resolution for every squarefree monomial ideal generated in degree $n-2$ in the polynomial ring $S=K[x_1,\ldots,x_n]$. Such an ideal is the complementary edge ideal $I_c(G)$ of a simple graph $G$ on $[n]$, and the resolution is supported on a labeled regular CW complex $X_G$ of dimension at most $2$. This yields closed formulas for the differentials and for the multigraded Betti numbers of $I_c(G)$ in terms of graph invariants. We then determine the $a$-invariant of $S/I_c(G)$ in terms of the vertex connectivity of $G$. Finally, when $S/I_c(G)$ is Cohen--Macaulay, we describe the canonical module and its Hilbert series, and recover the known characterizations of the Gorenstein and level cases. 
\end{abstract}
	
	\maketitle\vspace*{-1.4em}
	\section*{Introduction}

    Let $K$ be a field and $S=K[x_1,\ldots,x_n]$ be a polynomial ring. Describing the minimal free resolution of a monomial ideal of $S$ explicitly is a difficult problem, and for squarefree monomial ideals it has been solved only for special classes. One of the most useful tools for producing resolutions with a combinatorial flavour is the notion of a cellular resolution due to Bayer and Sturmfels \cite{bs1998}. One labels the vertices of a CW complex by the generators of the ideal, and labels every other cell by the least common multiple of the labels of its vertices. Then the cellular chain complex is a free resolution as soon as certain subcomplexes are acyclic. In this situation the Betti numbers can be read off from the cells of the complex; see also \cite{ms2005}.
In this paper we construct such a resolution for the class of complementary edge ideals. Let $n\geq 3$ be an integer. For a simple graph $G$ on the vertex set $[n]$ with at least one edge, the \emph{complementary edge ideal} of $G$ is defined as
\[
I_c(G)=\big(\mathbf x_{[n]\setminus e} : e\in E(G)\big)\subset S,
\]
where $\mathbf x_F=\prod_{i\in F}x_i$. Every squarefree monomial ideal generated in degree $n-2$ arises in this way. 
 Complementary edge ideals were introduced independently by Hibi, Qureshi and Saeedi Madani \cite{hqm2026} and by Ficarra and the author \cite{fm2025}. Since then, their algebraic properties and those of their powers have been studied from several viewpoints, among them linearity of powers \cite{fmStanley-Reisner2025},
 Cohen–Macaulayness, Gorensteinness, sequentially Cohen–Macaulayness, matroidality, levelness, linear resolutions, and the regularity of ordinary powers \cite{fm2025,hqm2026},  
 the licci property \cite{bcj2026,l2026}, associated primes \cite{f2026}, homological shifts \cite{lwz}, the equality of ordinary and symbolic powers \cite{rs2025}, their symbolic powers  with emphasis on regularity and depth  \cite{fmmReg2026,HK}, and their Rees algebras and symbolic Rees algebras  in \cite{fmRees2025} and \cite{fmm2026}, respectively.
 \smallskip

It was shown in \cite{hqm2026} that if $\beta_{i,j}(I_c(G))\neq0$, then $(i,j)\in\{(0,n-2),(1,n-1),(1,n),(2,n)\}$. In particular, $\pd\, I_c(G)\le 2$. Also for disconnected forests with no isolated vertices, the same work shows  that $\beta_{1,n}(I_c(G))\neq0$, and that $\beta_{1,n-1}(I_c(G))\neq0$ if and only if $G$ is not a disjoint union of edges. It is known by \cite[Corollary~4.7]{fm2025} that the graded Betti numbers of $I_c(G)^k$ do not depend on the field $K$ for any $k\ge1$.   Furthermore, in~\cite{jrs2026} it was shown that 
for a tree $T$, the ideal $I_c(T)$  admits a minimal Barile-Macchia resolution, and a combinatorial formula for the graded Betti numbers of $I_c(T)$ was given. 

Our main result establishes that for any graph $G$, the complementary edge ideal $I_c(G)$ admits a minimal cellular resolution supported on a regular CW complex. Let $G$ be a graph on $[n]$ with $n\geq 3$ and $m\geq 1$ edges, let $n_+$ be the number of non-isolated vertices of $G$, and let $c$ be the number of connected components of $G$ that are not isolated vertices. With this notation, we can state the theorem as follows.
\smallskip

 \noindent\textbf{Theorem A} (Theorem~\ref{thm:cellular-res} and Corollary~\ref{cor:explicit-resolution}).
\emph{The ideal $I_c(G)$ has a minimal cellular resolution supported on a labeled regular CW complex $X_G$ of dimension at most $2$. It has the form}
\begin{multline*}
	0 \longrightarrow S(-\mathbf x_{[n]})^{m-n_++c}
	\xrightarrow{\ d_2\ }
	S(-\mathbf x_{[n]})^{c-1} \oplus \bigoplus_{i\in V_+} S(-\mathbf x_{[n]\setminus\{i\}})^{\deg_G(i)-1}
	\xrightarrow{\ d_1\ }\\
	\bigoplus_{e \in E(G)} S(-\mathbf x_{[n]\setminus e})
	\xrightarrow{\ \varepsilon\ }
	I_c(G) \longrightarrow 0,
\end{multline*}
\emph{and all the differentials are given by explicit formulas.}
\smallskip 

The complex $X_G$ has a simple description. Its vertices are the edges of $G$, labeled by the generators of $I_c(G)$. For each non-isolated vertex $i$ of $G$ we fix an edge $a(i)$ containing $i$ and join $a(i)$ to every other edge containing $i$ by a $1$-cell with label $\mathbf x_{[n]\setminus\{i\}}$. The resulting graph $H$ has exactly $c$ connected components. We connect them by $c-1$ additional $1$-cells, and we attach a $2$-cell along the fundamental cycle of every edge of $H$ outside a chosen spanning forest of $H$. The additional $1$-cells and all $2$-cells carry the label $\mathbf x_{[n]}$. The number of $2$-cells is $m-n_++c$, the cycle rank of $G$. In particular, $X_G$ is a tree when $G$ is a forest. By construction, $X_G$ is a regular CW complex. By the Bayer--Sturmfels criterion, it suffices to show that the subcomplexes $X_{\preceq\mathbf b}$ are empty or acyclic. They turn out to be empty, single points or stars, hence trees, except when $\mathbf b$ has full support, in which case $X_{\preceq\mathbf b}=X_G$ is acyclic because its $2$-cells are attached along a basis of the cycle space of its $1$-skeleton. The choice of the edges $a(i)$ and of the spanning forest affects only the shape of $X_G$, not the Betti numbers.

Theorem~A has several consequences. The Betti table of $I_c(G)$ is determined by $m$, $n_+$ and $c$, and the nonzero graded Betti numbers are
\[
\beta_{0,n-2}=m,\qquad \beta_{1,n-1}=2m-n_+,\qquad \beta_{1,n}=c-1,\qquad \beta_{2,n}=m-n_++c.
\]
In particular, $I_c(G)$ has a linear resolution if and only if $c=1$. Moreover, $\pd\, I_c(G)$ equals $1$ if $G$ is a forest with at least two edges, $2$ if $G$ contains a cycle, and $0$ otherwise. The Hilbert series of $S/I_c(G)$ is given in Corollary~\ref{cor:hilbert}.

Our second main result concerns the $a$-invariant of $S/I_c(G)$. Applying Hochster's formula and Alexander duality, we reduce the computation of the relevant local cohomology modules to reduced cohomology groups of induced subgraphs of $G$. Here $\kappa(G)$ denotes the vertex connectivity of $G$.

 \noindent\textbf{Theorem B} (Theorem~\ref{thm:ainv}).
\emph{Let $G$ be a graph on $n\ge3$ vertices. Then}
\[
a(S/I_c(G))=\begin{cases}-\kappa(G), & \text{if $G$ has no isolated vertex and $G\neq K_n$},\\ 0,&\text{if } G=K_n,\\ -n_+,&\text{if $G$ has an isolated vertex}.\end{cases}
\]
\smallskip
In particular, for a graph without isolated vertices the $a$-invariant vanishes if and only if $G$ is disconnected or complete.

Finally, we turn to the Cohen--Macaulay case. Recall from \cite{fm2025,hqm2026} that,  $S/I_c(G)$ is Cohen--Macaulay if and only if $G$ is either a complete graph, or a forest with no isolated vertices or a graph with one edge. In these cases, the dual of the resolution of Theorem~A is a minimal free resolution of the canonical module $\omega_{S/I_c(G)}$. We use it to write down an explicit presentation of $\omega_{S/I_c(G)}$, to compute its Hilbert series and its multigraded Hilbert function, and to determine the Cohen--Macaulay type of $S/I_c(G)$. We also recover the known characterizations of the Gorenstein and level cases \cite{fm2025,hqm2026}. The paper ends with an example in which the canonical module is computed explicitly.

The paper is organized as follows. In Section~\ref{sec1} we fix notation and recall cellular resolutions following Bayer and Sturmfels, together with the Stanley--Reisner terminology. In Section~\ref{sec2} we construct the complex $X_G$ and prove Theorem~A, and we then derive the consequences for the Betti table, the projective dimension and the Hilbert series. In Section~\ref{three} we prove Theorem~B and give a description of the canonical module, its Hilbert series and multigraded Hilbert function.

\section{Preliminaries}\label{sec1}

Throughout, $K$ is a field, $S=K[x_1,\dots,x_n]$ is a standard graded polynomial ring, and $\m=(x_1,\dots,x_n)$ is its graded maximal ideal.
For $F\subseteq[n]$ we write $\mathbf x_F=\prod_{i\in F}x_i$, with $\mathbf x_\emptyset=1$, and $S(-\mathbf x_F)$ stands for $S(-\mathbf a)$, where $\mathbf a\in\{0,1\}^n$ is the indicator vector of $F$.

For a graph $G$ we write $V(G)=[n]$ for its vertex set, $E(G)$ for its edge set,  $Z$ for the set of isolated vertices of $G$, $m=|E(G)|$, and $\deg_G(i)$ for the degree of a vertex $i$. We set
\[
V_+=\{i\in V(G):\deg_G(i)>0\},\qquad n_+=|V_+|,
\]
and we let $c=c(G)$ denote the number of connected components of $G$ that are not isolated vertices. For $W\subseteq V(G)$, the \emph{induced subgraph} $G_W$ has vertex set $W$ and edge set $\{\{u,v\}\in E(G):u,v\in W\}$. A \emph{spanning tree} of a connected graph $G$ is a subgraph $T$ that is a tree with $V(T)=V(G)$. For an edge $q\in E(G)\setminus E(T)$, the unique cycle of $T+q$ is the \emph{fundamental cycle} of $q$ with respect to $T$. For a non-complete graph $G$, the \emph{vertex connectivity} of $G$ is
\[
\kappa(G)=\min\{|F|: F\subseteq V(G),\ G_{V(G)\setminus F}\text{ is disconnected}\}.
\]

Let $G$ be a graph on $[n]$ with $n\ge3$ and $E(G)\neq\emptyset$. The \emph{complementary edge ideal} of $G$ is defined as
\[
I_c(G)=\big(\mathbf x_{[n]\setminus e}: e\in E(G)\big)\subset S.
\]


Let $M$ be a finitely generated $\mathbb Z^n$-graded $S$-module with minimal multigraded free resolution
\[
0\to F_p\to\cdots\to F_1\to F_0\to M\to0,\qquad F_i=\bigoplus_{\mathbf a\in\mathbb Z^n}S(-\mathbf a)^{\beta_{i,\mathbf a}(M)} .
\]
The numbers $\beta_{i,\mathbf a}(M)$ are the \emph{multigraded Betti numbers} of $M$, and $\beta_{i,j}(M)=\sum_{|\mathbf a|=j}\beta_{i,\mathbf a}(M)$ are its $\mathbb Z$-graded Betti numbers. They are displayed in the \emph{Betti table}, where $\beta_{i,i+j}(M)$ is placed in column $i$ and row $j$. The \emph{projective dimension} $\pd\, M$ is the largest $i$ with $F_i\neq0$.
Let $d=\dim M$.
The value $a(M)=\max\{j : H^d_{\mathfrak m}(M)_j\neq0\}$ is called the {\em $a$-invariant} of $M$.

We recall the framework of cellular resolutions, following Bayer--Sturmfels \cite{bs1998}.
Let $X$ be a CW complex together with a labeling of its vertices (i.e.\ $0$-cells) by a chosen generating set of a monomial ideal $I\subset S$.  For any vertex $v$, let $\lambda(v)\in S$ be the monomial label of $v$. The label of an arbitrary face $F$ of $X$ is then defined as
	\[
	\lambda(F) := \operatorname{lcm}\big(\lambda(v) : v \text{ a vertex of } F\big).
	\]

Given such a labeled CW complex $X$, fix $K$-bases $\{[F]\}$ of the reduced cellular chain groups $C_{\bullet}(X;K)$. The \emph{cellular chain complex of $X$ supported by $\lambda$} is the complex $F_\bullet(X)$ of $\mathbb Z^n$-graded free $S$-modules
\[
F_i(X) = \bigoplus_{\dim F = i} S(-\lambda(F)),\qquad
\partial(F) = \sum_{F' \text{ facet of } F} [F:F']\,\frac{\lambda(F)}{\lambda(F')}\,F',
\]
where $[F:F']\in\{\pm1\}$ is the incidence number in the ordinary (unlabeled) cellular chain complex of $X$, and $\varepsilon(v)=\lambda(v)$ for each vertex $v$.

The labeled cell complex $(X,\lambda)$ \emph{supports a cellular resolution} of $I$, if the complex $F_\bullet(X)$ above, augmented by $\varepsilon: F_0(X)\to I$, is a free resolution of $I$.  

A \emph{simplicial complex on $[n]$} is a family $\Delta$ of subsets of $[n]$ closed under taking subsets. The ideal $I_\Delta=({\bf x}_F : F\subseteq[n],\ F\notin\Delta)$ is called the Stanley--Reisner ideal of $\Delta$, and $K[\Delta]=S/I_\Delta$ is the Stanley--Reisner ring of $\Delta$.
For $F\in\Delta$, the {\em link} of $F$ is defined as
$$\operatorname{lk}_\Delta F=\{H\in\Delta : H\cap F=\emptyset,\ H\cup F\in\Delta\}.$$
For a simplicial complex $\Delta$ on a vertex set $V$, the  {\em Alexander dual} of $\Delta$ is the complex $\Delta^\vee=\{H\subseteq V : V\setminus H\notin\Delta\}$.

\section{Minimal cellular resolutions of complementary edge ideals}\label{sec2}

In this section, we construct a minimal cellular resolution of $I_c(G)$ using a labeled CW complex of dimension at most two, and derive explicit formulas for its differentials and multigraded Betti numbers.
We will use the following criterion of Bayer--Sturmfels, which reduces the question of whether a CW complex $X$ supports a resolution to a topological condition on its subcomplexes.

\begin{Theorem}[Bayer--Sturmfels criterion, {\cite[Proposition 1.2]{bs1998}}]\label{thm:bs}
Let $(X,\lambda)$ be a labeled regular CW complex, and let $I\subset S$ be the monomial ideal generated by the labels of the vertices of $X$. For $\mathbf b\in\mathbb Z_{\ge0}^n$, let
\[
X_{\preceq \mathbf b} := \{F \in X : \lambda(F) \text{ divides } \mathbf x^{\mathbf b}\}.
\]
Then the cellular complex $F_\bullet(X)$, augmented by $\varepsilon$, is a free resolution of $I$ if and only if $X_{\preceq \mathbf b}$ is empty or $K$-acyclic  for every $\mathbf b\in\mathbb Z_{\ge0}^n$, that is $\widetilde H_i(X_{\preceq\mathbf b};K)=0$ for all $i$.
\end{Theorem}


\begin{Theorem}\label{thm:cellular-res}
				For any graph $G$ the ideal $I_c(G)$ admits a minimal cellular resolution supported on a labeled regular CW complex $X_G$ of dimension at most two.
			\end{Theorem}
			
			\begin{proof}
				Let $V(G)=[n]$, $V_+=\{i\in V(G):\, \deg(i)>0\}$, $m=|E(G)|$, $c=c(G)$, and $n_+=|V_+|$.
				For every $i\in V_+$, fix an edge of $G$ containing $i$ and denote it by $a(i)=\{i,t_i\}$. Notice that if $a(i)=a(j)$ for distinct vertices $i,j$, then necessarily $a(i)=a(j)=\{i,j\}\in E(G)$. We construct an auxiliary graph $H$ with vertex set $V(H)=E(G)$ and edge set
				\[
				E(H)=\{\{a(i),e\}:i\in V_+,\ e\in E(G),\ i\in e,\ e\neq a(i)\}.
				\]
				We denote the edge $\{a(i),e\}$ of $H$ by $q(i,e)$.
				
				We claim that $H$ has exactly $c$ connected components. Let $G_0$ be a connected component of $G$ containing at least one edge. Any two edges of $G_0$ can be joined by a sequence of pairwise incident edges. Two edges $e,f$ of $G_0$ that are incident at a vertex $i$, are either  equal to $a(i)$ or adjacent in $H$ to $a(i)$, so $e,f$ lie in the same component of $H$. Conversely, every edge of $H$ joins two edges of $G$ lying in the same component of $G$. Hence $H$ has $c$ connected components.
				
For distinct $u,v\in V(H)=E(G)$ with $\{u,v\}\in E(H)$, the edges $u,v$ of $G$ meet in a single vertex, which we denote by $h(u,v)$, and $a(h(u,v))\in\{u,v\}$. 
Set $P=\big\{(i,e)\ :\ i\in V_+,\ e\in E(G),\ i\in e,\ e\neq a(i)\big\}$. We show that the map
$\varphi:P\longrightarrow E(H)$,
defined as $\varphi((i,e))=\{a(i),e\}$
is a bijection. Clearly, $\varphi$ is well-defined and surjective.  	
Suppose $\varphi((i,e))=\varphi((i',e'))$, say both equal to $\{u,v\}\in E(H)$. Then $\{a(i),e\}=\{a(i'),e'\}=\{u,v\}$. 
	Since $(i,e)\in P$, we have $i\in e$. This gives
$i\in a(i)\cap e= u\cap v$.
	Likewise, $i'\in a(i')\cap e'=u\cap v$.
	Since $u$ and $v$ are distinct edges of the simple graph $G$, we have $|u\cap v|\le 1$. Thus
	$u\cap v=\{i\}=\{i'\}$,
	and hence $i=i'$.
	Also $\{a(i),e\}=\{a(i'),e'\}=\{a(i),e'\}$, and $e'\neq a(i')=a(i)$ by definition of $P$. Hence, $e'=e$,  proving injectivity of $\varphi$. Consequently
	\[ 
	|E(H)|=|P|=\sum_{i\in V_+}\big(\deg_G(i)-1\big)=2m-n_+.
	\]
				
Let $H_1,\dots,H_c$ be the connected components of $H$. Choose spanning trees $T_\ell\subseteq H_\ell$ and set $T=T_1\cup\cdots\cup T_c$. For $q\in E(H)\setminus E(T)$, let $C_q$ denote the unique cycle in $T+q$. Then $C_q\subseteq H_\ell$ for the $\ell$ with $q\in E(H_\ell)$, and
				\[
				|E(H)\setminus E(T)|=|E(H)|-\big(|V(H)|-c\big)=(2m-n_+)-(m-c)=m-n_++c.
				\]
				
We construct a labeled CW complex $X_G$ as follows.
The $0$-cells of $X_G$ are the vertices of $H$, i.e.\ the edges $e\in E(G)$, labeled $\lambda(e)=\mathbf x_{[n]\setminus e}$. Each $q(i,e)\in E(H)$ is a $1$-cell with endpoints $a(i)$ and $e$, labeled as
				\[
				\lambda(q(i,e))=\mathbf x_{[n]\setminus\{i\}}=\operatorname{lcm}\big(\lambda(a(i)),\lambda(e)\big).
				\]
				 We choose $p_\ell\in V(H_\ell)$ for each $1\leq\ell\leq c$, and add $c-1$ connector $1$-cells $b_2,\dots,b_c$, where $b_\ell$ joins $p_1$ to $p_\ell$, labeled $\lambda(b_\ell)=\mathbf x_{[n]}$. Finally, for each $q\in E(H)\setminus E(T)$ we attach a $2$-cell $\sigma_q$ along the unique cycle $C_q$ of $T+q$, labeled $\lambda(\sigma_q)=\mathbf x_{[n]}$. Notice that every edge of $C_q$ has label $\mathbf x_{[n]\setminus\{i\}}$ for some $i$, and $\operatorname{lcm}$ of the labels around the cycle is $\mathbf x_{[n]}$. We fix an orientation of every edge of $X_G^{(1)}$, as follows. Each edge $\{u,v\}\in E(H)$ is oriented from $u=a(h(u,v))$ to $v$, and each connector $b_\ell$ is oriented from $p_1$ to $p_\ell$.  This orientation determines the ordinary cellular boundary map $\partial_1:C_1(X_G^{(1)};K)\to C_0(X_G;K)$, as $\partial_1(u\to v)=v-u$.

				We show that $X_G$ is acyclic.
				Let $T'$ be the graph obtained from $T$ by adding the edges $b_2,\dots,b_c$. Since each $T_\ell$ spans $H_\ell$ and the connectors $b_\ell$ join the $c$ resulting trees into one tree without creating a cycle, $T'$ is a spanning tree of the $1$-skeleton $X_G^{(1)}$, and $E(X_G^{(1)})\setminus E(T')=E(H)\setminus E(T)$. For $q\in E(H_\ell)\setminus E(T)$, the unique $T'$-path between the endpoints of $q$ stays inside $T_\ell$, so $C_q$ is exactly the fundamental cycle of $q$ with respect to  $T'$ in $X_G^{(1)}$.
				For $q \in E(H)\setminus E(T)$ with $C_q=(e_1,\dots,e_s)$, we set $q_j=\{e_j,e_{j+1}\}$  and  $\eta_j=+1$ if the traversal $e_j\to e_{j+1}$ agrees with $q_j$'s fixed orientation, and $\eta_j=-1$ otherwise. Here $e_{s+1}=e_1$. Then $z_q=\sum_{j=1}^{s}\eta_j\,q_j \in Z_1(X_G^{(1)};K) = \ker \partial_1$.
				The set $\{z_q\}_{q\in E(H)\setminus E(T)}$ is a $K$-basis of $Z_1(X_G^{(1)};K)$, as each $z_q$ has coefficient $\pm1$ on $q$ and $0$ on every other non-tree edge, giving linear independence, and
			\[
			\dim_K Z_1(X_G^{(1)};K)=|E(X_G^{(1)})|-|V(X_G^{(1)})|+1=|E(H)\setminus E(T)|
			\]
			matches the number of $z_q$'s, so they also span $Z_1(X_G^{(1)};K)$.
                
				The boundary map $\partial_2:C_2(X_G;K)\to C_1(X_G^{(1)};K)$ sends the basis element $\sigma_q$ to $z_q$. So $\partial_2$ is injective with $\operatorname{im}\partial_2=\ker\partial_1$. Thus $\widetilde H_2(X_G;K)=0$ and $\widetilde H_1(X_G;K)=0$. Also since $X_G$ is connected, $\widetilde H_0(X_G;K)=0$ too. So $X_G$ is $K$-acyclic.
				
			To show that $X_G$ supports a cellular free resolution of $I_c(G)$, by Theorem~\ref{thm:bs}, it is enough to show that for any $\mathbf b\in\mathbb Z_{\ge0}^n$, the complex $X_{\preceq\mathbf b}=\{F\in X_G:\lambda(F)\mid \mathbf x^{\mathbf b}\}$ is acyclic. Since all labels are squarefree, $X_{\preceq \mathbf b}$ depends only on $A=\operatorname{supp}(\mathbf x^{\mathbf b})$. So we may take $\mathbf x^{\mathbf b}=\mathbf x_A$. If $|[n]\setminus A|\ge 3$, then no $0$-cell label divides $\mathbf x_A$, and so $X_{\preceq\mathbf b}=\emptyset$.
If $[n]\setminus A=\{i,j\}$, then  $X_{\preceq\mathbf b}$ is the point $\{i,j\}$ if $\{i,j\}\in E(G)$, and the empty set, otherwise.
If $[n]\setminus A=\{i\}$, then  
a cell $F$ satisfies $\lambda(F)\mid \mathbf x_A$ if and only if $i\notin\operatorname{supp}\lambda(F)$. If $i\notin V_+$, then no $0$-cell qualifies and $X_{\preceq\mathbf b}=\emptyset$. If $i\in V_+$, then a $0$-cell $e$ qualifies if and only if $i\in e$, and a $1$-cell $q(i',e)$ qualifies if and only if $i'=i$ (forcing it to have $a(i)$ as an endpoint). Hence $X_{\preceq\mathbf b}$ consists exactly of the vertex $a(i)$, the vertices $e\ne a(i)$ with $i\in e$, and the edges $q(i,e)$ joining each such $e$ to $a(i)$. This is the star centered at $a(i)$, which is a tree, hence $K$-acyclic.
Finally, if	$A=[n]$, then $X_{\preceq\mathbf b}=X_G$, which is $K$-acyclic as was shown above.
Thus $X_G$ supports a cellular free resolution of $I_c(G)$.
				
By construction of $X_G$ via lcm labels, label divisibility $\lambda(F')\mid\lambda(F)$ holds for every two faces $F'\subsetneq F$. Since connector $1$-cells $b_2,\dots,b_c$ never occur in a fundamental cycle $C_q\subseteq H_\ell$ and the other $1$-cells have label $\lambda(q(i,e))=\mathbf x_{[n]\setminus\{i\}}$, we have $\lambda(F')\ne\lambda(F)$ whenever $F'\subsetneq F$. 	So the divisibility is in fact strict, and hence the cellular free resolution supported by $X_G$ is minimal.
\end{proof}

Theorem~\ref{thm:cellular-res} gives the following explicit description of the multigraded minimal free resolution and of the Betti table of a complementary edge ideal, as follows. 			
\begin{Corollary}\label{cor:explicit-resolution}
With notation as in the proof of Theorem~\ref{thm:cellular-res}, the minimal multigraded free resolution of $I_c(G)$ supported on $X_G$ is of the form
\begin{multline*}
			0 \longrightarrow S(-\mathbf x_{[n]})^{m-n_++c}
\xrightarrow{\ d_2\ }
S(-\mathbf x_{[n]})^{c-1} \oplus \bigoplus_{i\in V_+} S(-\mathbf x_{[n]\setminus\{i\}})^{\deg_G(i)-1}
\xrightarrow{\ d_1\ }\\
\bigoplus_{e \in E(G)} S(-\mathbf x_{[n]\setminus e})
\xrightarrow{\ \varepsilon\ }
I_c(G) \longrightarrow 0.
		\end{multline*}
Using the orientation of $X_G^{(1)}$ fixed in the proof of Theorem~\ref{thm:cellular-res}, the augmentation and differentials are given explicitly by
$\varepsilon([e])=\mathbf x_{[n]\setminus e}$ for every $e\in E(G)$, 
\[
d_1(b_\ell)=\mathbf x_{p_\ell}\,[p_\ell]-\mathbf x_{p_1}\,[p_1], \qquad d_1(e\to f)={\bf x}_{f\setminus\{h\}}\,[f]-{\bf x}_{e\setminus\{h\}}\,[e],
\]
where $\{h\}=e\cap f$, and for any $q\in E(H)\setminus E(T)$, writing $C_q=(e_1,\dots,e_s)$, $q_j=\{e_j,e_{j+1}\}$, $\{h_j\}=e_j\cap e_{j+1}$, and $\eta_j$ the sign fixed in the proof of Theorem~\ref{thm:cellular-res}, 
\[
d_2(\sigma_q)=\sum_{j=1}^{s}\eta_j\, x_{h_j}\, q_j.
\]
In particular, the Betti table of $I_c(G)$ is of  the form
\[
\begin{array}{c|ccc}
 & 0 & 1 & 2 \\ \hline
0 & 0 & 0 & 0 \\
\vdots & \vdots & \vdots & \vdots \\
n-3 & 0 & 0 & 0 \\
n-2 & m & 2m-n_+ & m-n_++c \\
n-1 & 0 & c-1 & 0
\end{array}
\]
\end{Corollary}
\medskip

The Betti diagram shows, in particular, that $I_c(G)$ has a linear resolution if and only if $c=1$. More generally, this characterization extends to all powers of $I_c(G)$ by \cite[Corollary 3.2]{fmStanley-Reisner2025}, see also \cite[Theorem 2.2]{hqm2026}. Moreover, Corollary~\ref{cor:explicit-resolution} generalizes \cite[Theorem 3.3]{hqm2026} and implies \cite[Corollary 3.4]{hqm2026} and \cite[Proposition 4.9]{lwz}.

\begin{Corollary}\label{cor-pd}
For any graph $G$ with $m$ edges we have
$$\pd \, I_c(G) =\ \begin{cases}
			1,&\textit{if $G$ is a forest and $m\geq 2$,}\\
			2,&\textit{if $G$ contains a cycle,}\\
			0,&\textit{if $m=1$.}
		\end{cases}$$
\end{Corollary}

\begin{proof}
It follows from the Betti table of $I_c(G)$ that $\pd \, I_c(G)\leq 1$ if and only if the cycle rank $m-n_++c$ of $G$ is zero, that is $G$ is a forest. Thus $\pd \, I_c(G)=1$ if and only if  $G$ is a forest and either $c>1$ or $2m\neq n_+$, which is equivalent to $G$ being a forest with at least two edges.
\end{proof}


The multigraded minimal free resolution given in Corollary~\ref{cor:explicit-resolution} gives the  Hilbert series of $S/I_c(G)$ as follows. 

\begin{Corollary}\label{cor:hilbert}
The Hilbert series of $S/I_c(G)$ is
$$
H_{S/I_c(G)}(t)
=
\frac{
1-mt^{n-2}+(2m-n_+)t^{n-1}+(n_+-m-1)t^n
}{(1-t)^n}.
$$
\end{Corollary}



\section{The $a$-invariant and the canonical module} \label{three}

Let $I \subset S$ be a monomial ideal. Since $S/I$ carries a natural $\mathbb{N}^n$-graded structure inherited from $S$, each local cohomology module $H^i_{\mathfrak{m}}(S/I)$ inherits a $\mathbb{Z}^n$-graded structure. Let $H^i_{\mathfrak{m}}(S/I)_{\mathbf{a}}$ denote its graded component in multidegree $\mathbf{a} \in \mathbb{Z}^n$.
In order to compute the $a$-invariant of a complementary edge ideal, we use the following standard fact.

\begin{Lemma}
[Hochster's formula {\cite[Theorem~5.3.8]{BH}}]
\label{lem:hochster}
Let $\Delta$ be a simplicial complex on $[n]$, and let
${\bf a}\in\mathbb Z^n$.
Then $H^i_{\mathfrak m}(K[\Delta])_{\bf a}=0$ unless
$a_j\leq 0$ for all $j$ and
$F=\operatorname{supp}_-({\bf a})\in\Delta$.
In this case,
\[
\dim_K H^i_{\mathfrak m}(K[\Delta])_{\bf a}
=
\dim_K\widetilde H_{i-|F|-1}
(\operatorname{lk}_\Delta F;K).
\]
\end{Lemma}

The following result describes the $a$-invariant of a complementary edge ideal in terms of graph invariants.

\begin{Theorem}\label{thm:ainv}
Let $G$ be a graph on $n\geq 3$ vertices. Then

$$
a(S/I_c(G))=
\begin{cases}
-\kappa(G), & \text{if $G$ has no isolated vertices and }G\neq K_n,\\[2mm]
0, & \text{if }G=K_n,\\[2mm]
-n_+, & \text{if $G$ has an isolated vertex}.
\end{cases}
$$

In particular, $a(S/I_c(G))$ is independent of the field $K$. Moreover, if $G$ has no isolated vertices, then
$$
a(S/I_c(G))=0
\quad\Longleftrightarrow\quad
G\text{ is disconnected or }G=K_n.
$$
Otherwise, $a(S/I_c(G))=-\kappa(G)\le -1.$
\end{Theorem}

\begin{proof}
Let $\Gamma$ be the pure $1$-dimensional simplicial complex on $[n]$ whose facets are the edges of $G$. Set
$\Delta=\Gamma^\vee$. Then
\[
I_\Delta=({\bf x}_{[n]\setminus A} : A\in\Gamma)=({\bf x}_{[n]\setminus e} : e\in E(G))=I_c(G),
\]
and $S/I_c(G)=K[\Delta]$.
Set $d=\dim S/I_c(G)$, and  $h=\min\{|A| : A\subseteq[n],\ A\notin\Gamma\}$. Since $F\in\Delta$ if and only if $[n]\setminus F\notin\Gamma$, we have
$d=\max\{|F|:F\in\Delta\}=n-h$.

The component $H^d_{\mathfrak m}(K[\Delta])_j$ is the direct sum of the $H^d_{\mathfrak m}(K[\Delta])_{{\bf a}}$ over all ${\bf a}\in\mathbb Z^n$ with $\sum_i a_i=j$.
If $H^d_{\mathfrak m}(K[\Delta])_{{\bf a}}\neq0$, then by Lemma~\ref{lem:hochster} we have $a_i\le0$ for all $i$, $F=\operatorname{supp}_-({\bf a})\in\Delta$, $\widetilde H_{d-|F|-1}(\operatorname{lk}_\Delta F;K)\neq0$ and $\sum_ia_i\le-|F|$.
Conversely, if $F\in\Delta$ and $\widetilde H_{d-|F|-1}(\operatorname{lk}_\Delta F;K)\neq0$, then ${\bf a}=-\sum_{j\in F}e_j$ satisfies $H^d_{\mathfrak m}(K[\Delta])_{{\bf a}}\neq0$ and $\sum_ia_i=-|F|$.
So
\begin{equation}\label{eq:a-hochster}
a(S/I_c(G))=-\min\bigl\{|F| : F\in\Delta,\ \widetilde H_{d-|F|-1}(\operatorname{lk}_\Delta F;K)\neq0\bigr\}.
\end{equation}
Let $F\in\Delta$, put $W=[n]\setminus F$, and let $\Gamma_W:=\{A\in\Gamma : A\subseteq W\}$. For $D\subseteq[n]$ we have $D\in\operatorname{lk}_\Delta F$ if and only if $D\subseteq W$ and $[n]\setminus(D\cup F)=W\setminus D\notin\Gamma$, that is $W\setminus D\notin\Gamma_W$, or equivalently $D\in (\Gamma_W)^\vee$. Hence
\[
\operatorname{lk}_\Delta F=(\Gamma_W)^\vee.
\]
With $F$ and $W$ as above, put $i=d-|F|-1$. Since $|W|=n-|F|$ and $d=n-h$, we get $|W|-i-3=n-d-2=h-2$. Alexander duality  gives
\[
\widetilde H_{d-|F|-1}(\operatorname{lk}_\Delta F;K)\cong\widetilde H^{\,h-2}(\Gamma_W;K).
\]
Using this and \eqref{eq:a-hochster} we obtain
\begin{equation}\label{eq:a-final}
a(S/I_c(G))=-\min\bigl\{\,n-|W| : W\subseteq[n],\ W\notin\Gamma,\ \widetilde H^{\,h-2}(\Gamma_W;K)\neq0\bigr\}.
\end{equation}

{\bf Case 1}. Assume that $G$ has no isolated vertices and $G\neq K_n$. Then $h=2$, and  $\Gamma_W$ is the graph $G_W$ regarded as a simplicial complex. Hence $\widetilde H^0(\Gamma_W;K)\neq0$ if and only if $G_W$ is disconnected. Every such $W$ automatically satisfies $W\notin\Gamma$. Thus \eqref{eq:a-final} becomes
\[
a(S/I_c(G))=-\min\{|F| : G\setminus F\text{ is disconnected}\}=-\kappa(G).
\]
The set over which the minimum is taken is nonempty, as $G\neq K_n$ and there are non-adjacent vertices $u,v$, and so $F=[n]\setminus\{u,v\}$ qualifies.
\smallskip

{\bf Case 2}. Assume that $G=K_n$. Then  $h=3$, and $\Gamma_W$ is the complete graph on $W$, regarded as a $1$-dimensional complex. If $|W|\le2$, then $\Gamma_W$ is $\{\emptyset\}$, a point, or an edge. So $\widetilde H^1(\Gamma_W;K)=0$. If $|W|\ge3$, then $\Gamma_W$ is a connected graph with $\binom{|W|}{2}$ edges. So $\dim_K\widetilde H^1(\Gamma_W;K)=\binom{|W|}{2}-|W|+1=\binom{|W|-1}{2}>0$, and $W\notin\Gamma$. Thus the minimum in \eqref{eq:a-final} is attained at $W=[n]$. Hence $a(S/I)=0$.

\smallskip
{\bf Case 3}. Assume that $G$ has an isolated vertex. 
Let $Z$ be the set of isolated vertices of $G$, and $n_0=|Z|\ge1$. We have $h=1$, and  $\widetilde H^{-1}(\Gamma_W;K)\neq0$ if and only if $\Gamma_W=\{\emptyset\}$, i.e.\ if and only if $W\subseteq Z$. Moreover, for any $W\subseteq Z$, we have $W\notin\Gamma$ if and only if $W\neq\emptyset$. The minimum of $n-|W|$ over the nonempty subsets $W\subseteq Z$ is $n-n_0$, attained at $W=Z$. Hence $a(S/I)=-(n-n_0)=-n_+$.
\end{proof}

In the following we describe the canonical module $\omega_{S/I_c(G)}$. 

\begin{Remark}\label{rem:canonical}
{\em Let $G$ be a graph on $n\geq 3$ vertices, such that $R=S/I_c(G)$ is Cohen--Macaulay. It follows from \cite[Theorem 2.8]{fm2025} and \cite[Corollary 1.5]{hqm2026} that $G=K_n$, $G=K_2\sqcup (n-2)K_1$ or $G$ is a forest with no isolated vertices. Put $h=\operatorname{ht}I_c(G)=\pd_SR$ and $\omega_S=S(-\mathbf x_{[n]})$. The multigraded dual, with
 	respect to $\omega_S$, of the minimal free resolution of $R$
 	obtained from Corollary~\ref{cor:explicit-resolution} by
 	adjoining $S$ in homological degree $0$ is a minimal free
 	resolution of the canonical module
 	$\omega_R=\operatorname{Ext}^h_S(R,\omega_S)$. In particular,
 	$\omega_R$ is the cokernel of the last nonzero differential
 	in this dual complex. More precisely,
		\begin{itemize}
 		\item[\textup{(a)}] If $G$ is a forest with no isolated vertices,
 		then $h=2$ and 
 		$\omega_R=\operatorname{coker}(d_1^\vee)$, where
 		$d_1^\vee$ is the transpose of the matrix of $d_1$ in the
 		chosen homogeneous bases. The module $\omega_R$ has
 		$m-1=n-c-1$ minimal generators, corresponding to the
 		$1$-cells $q(i,e)$ and the connectors $b_\ell$ of $X_G$, and
 		its relations are given by the columns of $d_1^\vee$.
 		
 		\item[\textup{(b)}] If $G=K_n$, then $h=3$ and $c=1$. In this
 		case, $\omega_R=\operatorname{coker}(d_2^\vee)$, where
 		$d_2^\vee$ is the transpose of the matrix of $d_2$ in the
 		chosen homogeneous bases. The module $\omega_R$ has
 		$\binom{n-1}{2}$ minimal generators, corresponding to
 		the $2$-cells $\sigma_q$ of $X_G$.

        \item[\textup{(c)}] If $G=K_2\sqcup (n-2)K_1$, then $R=S/({\bf x}_Z)$ is a hypersurface ring, where $Z$ is the set of isolated vertices of $G$. Let $e$ be the unique edge of $G$. Then $\omega_R\cong (\omega_S/{\bf x}_Z\omega_S)({\bf 1}_Z)\cong (S/({\bf x}_Z))({\bf 1}_{Z}-{\bf 1}_{[n]})=R(-{\bf 1}_{e})$, see \cite[Theorem 3.54]{HMR}. 
 	\end{itemize}
	}
\end{Remark}

\begin{Proposition}
    \label{cor:canonical-invariants}
Assume that
$R=S/I_c(G)$ is Cohen--Macaulay. Then
\begin{enumerate}
\item[\textup{(a)}] The Cohen--Macaulay type of $R$ is 
\[
r(R)=\begin{cases}
m-1, & \text{if $G$ is a forest with no isolated vertex},\\[1mm]
\binom{n-1}{2}, & \text{if $G=K_n$},\\[1mm]
1, & \text{if $G=K_2\sqcup (n-2)K_1$}.
\end{cases}
\]
\item[\textup{(b)}] \textnormal{(}\cite[Theorem 2.12]{fm2025}, \cite[Corollary 1.6]{hqm2026}\textnormal{)} $R$ is a Gorenstein ring if and only if $G\in\{K_3,\,P_3,\,2K_2, K_2\sqcup (n-2)K_1\}$.
\item[\textup{(c)}] \textnormal{(}\cite[Corollary 3.5]{hqm2026}\textnormal{)} $R$ is level if and only if $G=K_n$, $G$ is a tree,  $G$ is a disjoint union of edges, or  $G=K_2\sqcup (n-2)K_1$.

	\item[(d)] The Hilbert series of $\omega_R$ is
 		\[
 		H_{\omega_R}(t)=
 		\begin{cases}
 			\dfrac{(c-1)+\sum_{j=1}^{n-2}(n-1-j)\,t^{j}}{(1-t)^{n-2}}, & \text{if $G$ is a forest with no isolated vertices},\\[4mm]
 			\dfrac{\sum_{j=0}^{n-3}\binom{n-1-j}{2}\,t^{j}}{(1-t)^{n-3}}, & \text{if } G=K_n,\\[4mm]
         \dfrac{t^2(1+t+\cdots+t^{n-3})}{(1-t)^{n-1}},  & \text{if } G=K_2\sqcup (n-2)K_1.
 		\end{cases}
 		\]
		    
\item[\textup{(e)}] $(\omega_R)_{{\bf b}}=0$ for ${\bf b}\notin\mathbb Z_{\ge0}^n$. For ${\bf b}\in\mathbb Z_{\ge0}^n$, put $F=\operatorname{supp}({\bf b})$
and $W=[n]\setminus F$. If $G$ is a forest with no isolated vertices, then
\[
\dim_K(\omega_R)_{{\bf b}}=\begin{cases}\operatorname{comp}(G_W)-1,& W\neq\emptyset,\\ 0,& W=\emptyset,\end{cases}
\]
where $\operatorname{comp}(G_W)$ denotes the number of connected components of $G_W$. If $G=K_n$, then
\[
\dim_K(\omega_R)_{{\bf b}}=\begin{cases}\binom{n-1-|F|}{2},& |F|\le n-1,\\ 0,& |F|=n.\end{cases}
\]
\end{enumerate}

If $G=K_2\sqcup (n-2)K_1$, then 
\[
\dim_K(\omega_R)_{{\bf b}}=\begin{cases}1,& \emptyset\neq W\subseteq Z,\, \\ 0,& \text{otherwise}.\end{cases}
\]
\end{Proposition}

\begin{proof}
(a) Note that $r(R)=\mu(\omega_R)$ is the rank of the last module of the resolution  given in Corollary~\ref{cor:explicit-resolution}, namely
$(c-1)+\sum_i(\deg_G(i)-1)=2m-n+c-1=m-1$ for a forest, and $m-n+1=\binom{n-1}{2}$ for $K_n$. If $G=K_2\sqcup (n-2)K_1$, then  $R$ is a hypersurface. So it has type 1, and it is Gorenstein and level.
\smallskip

(b) $R$ is Gorenstein if and only if $r(R)=1$. Now $\binom{n-1}{2}=1$ if and only if $n=3$, and $m-1=1$ if and only if $m=2$. A forest with two edges
and no isolated vertices is $P_3$ or $2K_2$.
\smallskip

(c) $R$ is level if and only if  $\omega_R$ is generated in a single degree. For $K_n$ this holds, by Remark~\ref{rem:canonical}. For a forest, $\omega_R$ has
$c-1$ generators in degree $0$ and $\sum_i(\deg_G(i)-1)=n-2c$ generators in degree $1$. So $R$ is level if and only if $c=1$
or $n=2c$. This is equivalent to $G$ being a tree or every component of $G$ is an edge.
\smallskip

(d) Let $G$ be a forest with no isolated vertices. Then $2m-n=n-2c$ and $n=n_+$. Corollary~\ref{cor:explicit-resolution} and Remark~\ref{rem:canonical} yield that
$H_{\omega_R}(t)=g(t)/(1-t)^n$, where 
$g(t)=(c-1)+(n-2c)\,t-(n-c)\,t^{2}+t^{n}$.
	Put $q_0=c-1$, $q_j=n-1-j$ for $1\le j\le n-1$ (so $q_{n-1}=0$), and $q_j=0$ for $j<0$ or $j\ge n$. Then $Q(t)=\sum_{j\ge0}q_jt^j=(c-1)+\sum_{j=1}^{n-2}(n-1-j)t^j$. We show that  $(1-t)^2Q(t)=g(t)$. The coefficient of $t^k$ in $(1-t)^2Q(t)$ is $s_k=q_k-2q_{k-1}+q_{k-2}$, and
	\begin{itemize}
		\item $s_0=c-1$ and $s_1=(n-2)-2(c-1)=n-2c$,
		\item $s_2=(n-3)-2(n-2)+(c-1)=-(n-c)$,
		\item for $3\le k\le n-1$, the indices $k,k-1,k-2$ lie in $[1,n-1]$, where $q_j$ is linear in $j$, so $s_k=0$,
		\item $s_n=0-0+q_{n-2}=1$, and $s_k=0$ for $k>n$.
	\end{itemize}
	Hence $(1-t)^2Q(t)=g(t)$. Therefore
	\[
	H_{\omega_R}(t)=\frac{g(t)}{(1-t)^n}=\frac{Q(t)}{(1-t)^{n-2}}=\frac{(c-1)+\sum_{j=1}^{n-2}(n-1-j)t^j}{(1-t)^{n-2}}.
	\]
Now, assume that $G=K_n$. Then $H_{\omega_R}(t)=f(t)/(1-t)^n$, where 
$f(t)=\tbinom{n-1}{2}-n(n-2)\,t+\tbinom n2\,t^{2}-t^{n}$.	
	Put $r_j=\binom{n-1-j}{2}$ for $0\le j\le n-1$, and $r_j=0$ otherwise. Note that $r_{n-2}=r_{n-1}=0$. Set $Q(t)=\sum_{j\ge0}r_jt^j=\sum_{j=0}^{n-3}\binom{n-1-j}{2}t^j$. We show that $(1-t)^3Q(t)=f(t)$. The coefficient of $t^k$ in $(1-t)^3Q(t)$ is $s_k=r_k-3r_{k-1}+3r_{k-2}-r_{k-3}$, and 
	\begin{itemize}
		\item $s_0=\binom{n-1}{2}$,
		\item $s_1=\binom{n-2}{2}-3\binom{n-1}{2}=-n(n-2)$,
		\item $s_2=\tfrac12\big[(n-3)(n-4)-3(n-2)(n-3)+3(n-1)(n-2)\big]=\tfrac12 n(n-1)=\binom n2$,
		\item for $3\le k\le n-1$, all indices $k,\dots,k-3$ lie in $[0,n-1]$, where $r_j$ is a quadratic polynomial in $j$, and the third finite difference of a quadratic vanishes, so $s_k=0$,
		\item $s_n=0-3r_{n-1}+3r_{n-2}-r_{n-3}=-\binom22=-1$, and $s_k=0$ for $k>n$.
	\end{itemize}
	Hence $(1-t)^3Q(t)=f(t)$. Therefore
	\[
	H_{\omega_R}(t)=\frac{f(t)}{(1-t)^n}=\frac{Q(t)}{(1-t)^{n-3}}=\frac{\sum_{j=0}^{n-3}\binom{n-1-j}{2}t^j}{(1-t)^{n-3}}. 
	\]
    Finally, let $G=K_2\sqcup (n-2)K_1$. By Remark~\ref{rem:canonical}(c), we have $\omega_R\cong R(-{\bf 1}_{e})$, which gives $H_{\omega_R}(t)=t^2H_{R}(t)=\frac{t^2(1-t^{n-2})}{(1-t)^n}=\frac{t^2(1+t+\cdots+t^{n-3})}{(1-t)^{n-1}}.$ 
\smallskip

(e) Let $F_\bullet$ be the free resolution given in Corollary \ref{cor:explicit-resolution}. Clearly, all generators of $F_\bullet^\vee$ have degrees in $\{0,1\}^n$, and so $(\omega_R)_{{\bf b}}=0$ for ${\bf b}\notin\mathbb Z_{\ge0}^n$.
For ${\bf b}\in\mathbb Z_{\ge0}^n$ and squarefree $\mathbf x_A$ we have $\dim_KS(-\mathbf x_A)_{{\bf b}}=1$ if $A\subseteq F$ and $0$ otherwise,
and $\dim_K(\omega_R)_{{\bf b}}$ is the alternating sum of the dimensions in $F_\bullet^\vee$. For $B\subseteq V(G)$, we set $e(B)=|E(G_B)|$, and for disjoint subsets $B,C\subseteq V(G)$, let
$e(B,C)=|\{\{u,v\}\in E(G):u\in B,\ v\in C\}|$.
Then for a forest $G$, we have
\[
\dim_K(\omega_R)_{{\bf b}}=(c-1)+\sum_{i\in F}(\deg_G(i)-1)-e(F)+\delta_{F,[n]} .
\]
Since $\sum_{i\in F}\deg_G(i)=2e(F)+e(F,W)$, $m=e(F)+e(F,W)+e(W)$ and $m=n-c$, we obtain 
$\dim_K(\omega_R)_{{\bf b}}=|W|-e(W)-1+\delta_{W,\emptyset}$. As $G_W$ is a forest, $|W|-e(W)=\operatorname{comp}(G_W)$, which gives the claim.

For $K_n$ one gets
$\dim_K(\omega_R)_{{\bf b}}=\binom{n-1}{2}-(n-2)|F|+\binom {|F|}{2}-\delta_{F,[n]}$, which equals $\binom{n-1-|F|}{2}$ for
$|F|\le n-1$ and $0$ for $F=[n]$.

For $G=K_2\sqcup (n-2)K_1$, we have $\dim_K(\omega_R)_{{\bf b}}=\dim S_{\mathbf b-\mathbf 1_e}-\dim S_{\mathbf b-\mathbf 1_{[n]}}$.
If $\mathbf b-\mathbf 1_e\notin\mathbb Z_{\ge0}^n$, both vector space dimensions are $0$ and $(\omega_R)_{\mathbf b}=0$.
Now, let $\mathbf b-\mathbf 1_e\in\mathbb Z_{\ge0}^n$. Then $\dim(\omega_R)_{\mathbf b}=0$ if $b_z\ge1$ for every $z\in Z$, and $\dim(\omega_R)_{\mathbf b}=1$ if $b_z=0$ for some $z\in Z$. Note that since ${\bf b}\in\mathbb Z_{\ge0}^n$, we have $b_z=0$ if and only if $z\in W$. Also  $\mathbf b-\mathbf 1_e\in\mathbb Z_{\ge0}^n$ if and only if $W\subseteq Z$, which completes the proof.
\end{proof}

\begin{Example}
{\em Let $G=P_4$ with edges $f_1=\{1,2\}$, $f_2=\{2,3\}$, $f_3=\{3,4\}$, let $S=K[x_1,\dots,x_4]$, and $R=S/I_c(G)$, so that $I_c(G)=(x_3x_4,\,x_1x_4,\,x_1x_2)$. Choose $a(1)=a(2)=f_1$, $a(3)=f_2$ and $a(4)=f_3$. Then the $1$-cells of $H$ are $f_1\to f_2$ (with $h(f_1,f_2)=2$) and $f_2\to f_3$ (with $h(f_2,f_3)=3$). Put $\alpha=(f_1\to f_2)^*$ and $\beta=(f_2\to f_3)^*$. By the formulas of Corollary~\ref{cor:explicit-resolution},
\[
d_1(f_1\to f_2)=x_3[f_2]-x_1[f_1],\qquad d_1(f_2\to f_3)=x_4[f_3]-x_2[f_2],
\]
and transposing gives
\[
\omega_R=\operatorname{coker}(d_1^*)\;\cong\;(S\alpha\oplus S\beta)\big/\big(x_1\alpha,\ x_3\alpha-x_2\beta,\ x_4\beta\big),
\qquad \deg\alpha=\mathbf e_2,\ \deg\beta=\mathbf e_3,
\]
where $\mathbf e_i\in\mathbb Z^4$ denotes the $i$-th standard basis vector. The three relations have multidegrees $\mathbf 1_{f_1}$, $\mathbf 1_{f_2}$, $\mathbf 1_{f_3}$, respectively. In the standard $\mathbb Z$-grading, $\deg\alpha=\deg\beta=1$, and
\[
H_{\omega_R}(t)=\frac{2t+t^{2}}{(1-t)^{2}}.
\]
In particular, $r(R)=m-1=2$ and $a(R)=-1=-\kappa(P_4)$. So $R$ is Cohen--Macaulay but not Gorenstein.
}
\end{Example}

In  \cite{fm2025}  it was shown that for any $k\geq 1$,   $\operatorname{reg} I_c(G)^k=(n-2)k+\min\{k,c-1\}$  and that the
graded Betti numbers of $I_c(G)^k$ do not depend on $K$. Moreover,
$I_c(G)^k$ has linear quotients for all $k\ge1$ if and only if
$I_c(G)$ has a linear resolution, i.e.\ $c=1$, see \cite[Corollary 3.2]{fmStanley-Reisner2025}. 
So for connected $G$ the Betti numbers of $I_c(G)^k$ can be computed from linear quotients. In light of these facts and Theorem~\ref{thm:cellular-res}  we conclude the paper with the following question.

\begin{Question}
Let $G$ be a graph on $[n]$ and $k\ge2$. Does $I_c(G)^k$ admit a minimal cellular
resolution? 
\end{Question}

			\bigskip
			
\textbf{Acknowledgment.} This work is based upon research funded by Iran National Science Foundation (INSF) under project No. 40507515.

		\end{document}